\documentclass[12pt, reqno]{amsart}
\usepackage{amsmath, amsthm, amscd, amsfonts, amssymb, graphicx, color}
\usepackage[bookmarksnumbered, colorlinks, plainpages]{hyperref}

\newtheorem{theorem}{Theorem}[section]
\newtheorem{lemma}[theorem]{Lemma}

\theoremstyle{definition}

\theoremstyle{remark}

\numberwithin{equation}{section}

\begin{document}
\setcounter{page}{1}

\title[Local derivations on matrix algebras]{Local derivations on
associative and Jordan matrix algebras}

\author[Sh.\ Ayupov]{Shavkat Ayupov$^{1,2}$}

\address{$^{1}$
V.I. Romanovskiy Institute of Mathematics, Uzbekistan Academy of Sciences,
Tashkent, Uzbekistan.}
\address{$^{2}$
National University of Uzbekistan,
Tashkent, Uzbekistan.}
\email{\textcolor[rgb]{0.00,0.00,0.84}{sh$_-$ayupov@mail.ru}}

\author[F. Arzikulov]{Farhodjon Arzikulov$^3$}

\address{$^{3}$ Department of Mathematics, Andizhan State University, Andizhan, Uzbekistan.}
\email{\textcolor[rgb]{0.00,0.00,0.84}{arzikulovfn@rambler.ru}}

%\dedicatory{This paper is dedicated to Professor ABCD}



\thanks{}

\begin{abstract}
In the present paper we prove that
every additive (not necessarily homogenous) local inner
derivation on the algebra of matrices over an arbitrary
field is an inner derivation, and every local inner
derivation on the ring of matrices over a finite ring generated by the identity element
 or
the ring of integers is an inner derivation.
We also prove that every additive
local inner derivation on the Jordan algebra of symmetric matrices over an arbitrary
field is a derivation, and every local inner
derivation on the Jordan ring of symmetric matrices over a finite ring generated by the identity element
 or
the ring of integers is a derivation.
\end{abstract}

\maketitle

\section{Introduction}

The present paper is devoted to local derivations on associative
and Jordan matrix algebras. Recall that a local derivation is defined as follows:
given an algebra $\mathcal{A}$, a linear map $\nabla : \mathcal{A} \to \mathcal{A}$
is called a local derivation if for every $x\in
\mathcal{A}$ there exists a derivation $D : \mathcal{A}\to \mathcal{A}$ such that
$\nabla(x)=D(x)$.

In \cite{1}, R. Kadison introduces the concept of local derivation and proves
that each continuous local derivation from a von Neumann algebra into its dual Banach
bemodule is a derivation. B. Jonson \cite{2} extends the above result by proving that every
local derivation from a C*-algebra into its Banach bimodule is a derivation. In particular, Johnson
gives an automatic continuity result by proving that local derivations of a C*-algebra $A$ into a
Banach $A$-bimodule $X$ are continuous even if not assumed a priori to be so
(cf. \cite[Theorem 7.5]{2}). Based on these results, many authors have studied
local derivations on operator algebras, for example, see in 
\cite{AA4, SAK, AK3, AK4, CLW, 3, HL1, 4, HLAH, HLZ, J, KK, 5, LP, LZ1, LZ2, 6, ZJC, ZPY}.

In this paper we develop a pure algebraic approach to
investigation of derivations and local derivations on
associative and Jordan algebras. Since we consider a sufficiently general case
we restrict our attention only on inner
derivations and local inner derivations.

In section 2 we introduce and investigate a notion of additive local derivation on the
algebra $M_n(\mathcal{F})$ of matrices over an arbitrary field $\mathcal{F}$. It is proved that, given
an arbitrary field $\mathcal{F}$, every additive local inner
derivation on the algebra $M_n(\mathcal{F})$ is an inner derivation, and every local inner
derivation on the ring of matrices over a finite ring generated by the identity element
 or
the ring of integers is an inner derivation, where a finite ring is a ring that has a finite number of elements.
Here we define an additive local inner
derivation as follows: given an algebra $\mathcal{A}$, an additive (not necessarily homogenous) map $\nabla : \mathcal{A} \to \mathcal{A}$
is called additive local inner derivation if for every $x\in
\mathcal{A}$ there exists an inner derivation $D : \mathcal{A}\to \mathcal{A}$ such that
$\nabla(x)=D(x)$.

In section 3 additive local derivations on the Jordan algebra of symmetric matrices over an arbitrary field
are introduced and studied. It is proved that every additive local inner derivation on the Jordan
algebra $H_n(\mathcal{F})$ of $n$-dimensional symmetric matrices over an arbitrary field $\mathcal{F}$ is
a derivation, and every local inner
derivation on the Jordan ring of symmetric matrices over a finite ring generated by the identity element
 or
the ring of integers is a derivation. For this propose we use a Jordan analogue of the algebraic approach to the
investigation of additive local derivations applied to
algebras of matrices over an arbitrary field developed in section 2.
The method developed in this paper is sufficiently universal
and can also be applied to Jordan and Lie algebras. Its corresponding
modification has been used when we considered a similar problem for Lie algebras of
skew-symmetric matrices over an arbitrary field \cite{AA4}.
It should be noted that the notions of local inner derivation and local spatial derivation in theorems
5.4, 5.5 and 6.1 in \cite{AA4} can be replaced by the
notions of additive local inner derivation and additive local spatial derivation respectively.

\section{local derivations on associative algebras of matrices}

Let $\mathcal{A}$ be an algebra. A liner map $D : \mathcal{A}\to \mathcal{A}$ is called
a derivation, if $D(xy)=D(x)y+xD(y)$ for any two elements $x$, $y\in \mathcal{A}$.

An additive map $\nabla : \mathcal{A}\to \mathcal{A}$ is called additive local derivation, if for
any element $x\in \mathcal{A}$ there exists a derivation
$D:\mathcal{A}\to \mathcal{A}$ such that $\nabla (x)=D(x)$.

Now let $\mathcal{A}$ be a non commutative (but associative)  algebra. A derivation $D$ on $\mathcal{A}$
is called an inner derivation, if there exists an element $a\in \mathcal{A}$ such that
$$
D(x)=ax-xa, x\in \mathcal{A}.
$$
This derivation $D$ we denote by $D_a$, i.e. $D_a(x)=ax-xa$.
An additive (not necessarily homogenous) map $\nabla : \mathcal{A}\to \mathcal{A}$ is called additive local inner derivation,
if for any element $x\in \mathcal{A}$ there exists an inner derivation $D_a$ such that
$\nabla(x)=D_a(x)$.

Let $\mathcal{F}$ be a field, and let $M_n(\mathcal{F})$ be the
matrix algebra over $\mathcal{F}$, $n>1$, i.e. consisting of matrices
$$
\left[%
\begin{array}{cccc}
 a^{1,1} & a^{1,2} & \cdots & a^{1,n}\\
a^{2,1} & a^{2,2} & \cdots & a^{2,n}\\
\vdots & \vdots & \ddots & \vdots\\
a^{n,1} & a^{n,2} & \cdots & a^{n,n}\\
\end{array}%
\right],
a^{i,j}\in \mathcal{F}, i,j=1,2,\dots ,n.
$$
Let $\{e_{i,j}\}_{i,j=1}^n$ be the set of matrix units in
$M_n(\mathcal{F})$, i.e. $e_{i,j}$ is the matrix with components
$a^{i,j}={\bf 1}$ and $a^{k,l}={\bf 0}$ if $(i,j)\neq(k,l)$, where
${\bf 1}$ is the identity element, ${\bf 0}$ is the zero element of the field
$\mathcal{F}$, and a matrix $a\in M_n(\mathcal{F})$ is written as $a=\sum_{k,l=1}^n
a^{k,l}e_{k,l}$, where $a^{k,l}\in \mathcal{F}$ for $k,l=1,2,\dots, n$,
or as $a=\sum_{k,l=1}^n a_{k,l}$, where $a_{k,l}=e_{k,k}ae_{l,l}$ for $k,l=1,2,\dots, n$.

\medskip

First, let us prove some lemmas which will be used in the
proof of Theorem \ref{14}. Throughout the section, $\mathcal{F}$ denotes an arbitrary field, $M_n(\mathcal{F})$ denotes the algebra of $n\times
n$ matrices over $\mathcal{F}$, $n>1$. Let  $\nabla :M_n(\mathcal{F})\to M_n(\mathcal{F})$ be an additive local inner derivation.

\medskip

\begin{lemma} \label{1}
For arbitrary $\lambda$, $\mu\in \mathcal{F}$ and each pair $i$, $j$ of distinct indices there exists an element $a\in M_n(\mathcal{F})$ such that
$$
\nabla(\lambda e_{i,i})=D_a(\lambda e_{i,i}), \nabla(\mu e_{i,j})=D_a(\mu e_{i,j}).
$$
\end{lemma}

\begin{proof}
We have
$$
\nabla(\lambda e_{i,i})+\nabla(\mu e_{i,j})=\nabla(\lambda e_{i,i}+\mu e_{i,j}).
$$
Let $d_{i,i}$, $d_{i,j}$, $d\in M_n(\mathcal{F})$ be elements such that
$$
\nabla(\lambda e_{i,i})=d_{i,i}\lambda e_{i,i}-\lambda e_{i,i}d_{i,i},
\nabla(\mu e_{i,j})=d_{i,j}\mu e_{i,j}-\mu e_{i,j}d_{i,j},
$$
$$
\nabla(\lambda e_{i,i}+\mu e_{i,j})=d(\lambda e_{i,i}+\mu e_{i,j})-(\lambda e_{i,i}+\mu e_{i,j})d.
$$
Note that $e_{i,i}d_{i,j}e_{j,j}\mu e_{i,j}=\mu e_{i,j}e_{i,i}d_{i,j}e_{j,j}=0$. So we can take
$e_{i,i}d_{i,j}e_{j,j}=e_{i,i}d_{i,i}e_{j,j}$, i.e. $d_{i,j}^{i,j}=d_{i,i}^{i,j}$.
Since
$$
e_{i,i}[d_{i,i}\lambda e_{i,i}-\lambda e_{i,i}d_{i,i}+d_{i,j}\mu e_{i,j}-\mu e_{i,j}d_{i,j}]e_{i,i}=e_{i,i}[d(\lambda e_{i,i}+\mu e_{i,j})-(\lambda e_{i,i}+\mu e_{i,j})d]e_{i,i}
$$
we have
$$
-\mu e_{i,j}d_{i,j}e_{i,i}=-\mu e_{i,j}de_{i,i},\,\,i.e \,\,d_{i,j}^{j,i}=d^{j,i}.
$$
The equalities
$$
e_{j,j}[d_{i,i}\lambda e_{i,i}-\lambda e_{i,i}d_{i,i}+d_{i,j}\mu e_{i,j}-\mu e_{i,j}d_{i,i}]e_{i,i}=e_{j,j}[d(\lambda e_{i,i}+\mu e_{i,j})-(\lambda e_{i,i}+\mu e_{i,j})d]e_{i,i}
$$
imply that
$$
e_{j,j}d_{i,i}e_{i,i}=e_{j,j}de_{i,i},\,\,i.e \,\,d_{i,i}^{j,i}=d^{j,i}.
$$
Hence
$$
d_{i,j}^{j,i}=d_{i,i}^{j,i}.
$$

Let $e=e_{i,i}+e_{j,j}$. Then
$$
(1-e)(D_d(\lambda e_{i,i}+\mu e_{i,j}))e_{i,i}=(1-e)(D_{d_{i,i}}(\lambda e_{i,i}))e_{i,i}+(1-e)(D_{d_{i,j}}(\mu e_{i,j}))e_{i,i}
$$
and
$$
(1-e)de_{i,i}=(1-e)d_{i,i}e_{i,i}.
$$
Also
$$
(1-e)(D_d(\lambda e_{i,i}+\mu e_{i,j}))e_{j,j}=(1-e)(D_{d_{i,i}}(\lambda e_{i,i}))e_{j,j}+(1-e)(D_{d_{i,j}}(\mu e_{i,j}))e_{j,j}
$$
and
$$
(1-e)de_{i,j}=(1-e)d_{i,j}e_{i,j},\,\, i.e.\,\, (1-e)de_{i,i}=(1-e)d_{i,j}e_{i,i}.
$$
Hence $(1-e)d_{i,i}e_{i,i}=(1-e)d_{i,j}e_{i,i}$.
At the same time, since $e_{i,i}d_{i,j}(1-e)e_{i,j}=e_{i,j}e_{i,i}d_{i,j}(1-e)=0$,
we can take $e_{i,i}d_{i,j}(1-e)=e_{i,i}d_{i,i}(1-e)$.

Therefore
$$
\nabla(\lambda e_{i,i})=e_{i,i}d_{i,i}\lambda e_{i,i}+e_{j,j}d_{i,i}\lambda e_{i,i}+(1-e)d_{i,i}\lambda e_{i,i}-\lambda e_{i,i}d_{i,i}e_{j,j}-
$$
$$
\lambda e_{i,i}d_{i,i}e_{i,i}-\lambda e_{i,i}d_{i,i}(1-e)=
$$
$$
e_{j,j}d_{i,j}\lambda e_{i,i}-\lambda e_{i,i}d_{i,j}e_{j,j}+e_{i,i}d_{i,j}\lambda e_{i,i}+(1-e)d_{i,j}\lambda e_{i,i}-
$$
$$
\lambda e_{i,i}d_{i,j}e_{i,i}-\lambda e_{i,i}d_{i,j}(1-e)=
$$
$$
d_{i,j}\lambda e_{i,i}-\lambda e_{i,i}d_{i,j}=D_{d_{i,j}}(\lambda e_{i,i}).
$$
since
$$
e_{i,i}d_{i,i}\lambda e_{i,i}-\lambda e_{i,i}d_{i,i}e_{i,i}=e_{i,i}d_{i,j}\lambda e_{i,i}-\lambda e_{i,i}d_{i,j}e_{i,i}=0,
$$
$$
d_{i,j}^{i,j}=d_{i,i}^{i,j}.
$$
This completes the proof.
\end{proof}

\medskip
Similarly we can prove the following lemma.

\begin{lemma} \label{2}

For arbitrary $\lambda$, $\mu\in \mathcal{F}$ and each pair $i$, $j$ of distinct indices there exists an element $a\in M_n(\mathcal{F})$ such that
$$
\nabla(\lambda e_{i,i})=D_a(\lambda e_{i,i}), \nabla(\mu e_{j,i})=D_a(\mu e_{j,i}).
$$
\end{lemma}

\begin{lemma} \label{3}
For arbitrary $\lambda$, $\mu\in \mathcal{F}$ and each pair $i$, $j$ of distinct indices there exists an element $a\in M_n(\mathcal{F})$ such that
$$
\nabla(\lambda e_{i,j})=D_a(\lambda e_{i,j}), \nabla(\mu e_{j,i})=D_a(\mu e_{j,i}).
$$
\end{lemma}

\begin{proof}
We have
$$
\nabla(\lambda e_{i,j})+\nabla(\mu e_{j,i})=\nabla(\lambda e_{i,j}+\mu e_{j,i})
$$
Let $d_{i,j}$, $d_{j,i}$, $d\in M_n(\mathcal{F})$ be such elements that
$$
\nabla(\lambda e_{i,j})=d_{i,j}\lambda e_{i,j}-\lambda e_{i,j}d_{i,j},
\nabla(\mu e_{j,i})=d_{j,i}\mu e_{j,i}-\mu e_{j,i}d_{j,i},
$$
$$
\nabla(\lambda e_{i,j}+\mu e_{j,i})=d(\lambda e_{i,j}+\mu e_{j,i})-(\lambda e_{i,j}+\mu e_{j,i})d.
$$
Then
$$
d_{i,j}\lambda e_{i,j}-\lambda e_{i,j}d_{i,j}+d_{j,i}\mu e_{j,i}-\mu e_{j,i}d_{j,i}=d(\lambda e_{i,j}+\mu e_{j,i})-(\lambda e_{i,j}+\mu e_{j,i})d.
$$
Note that
$$
e_{i,i}d_{i,j}e_{j,j}\lambda e_{i,j}=\lambda e_{i,j}e_{i,i}d_{i,j}e_{j,j}=e_{j,j}d_{j,i}e_{i,i}\mu e_{j,i}=\mu e_{j,i}e_{j,j}d_{j,i}e_{i,i}=0.
$$
So we can take
$$
e_{i,i}d_{i,j}e_{j,j}=e_{i,i}d_{j,i}e_{j,j},\,\, i.e.\,\, d_{i,j}^{i,j}=d_{j,i}^{i,j},
$$
$$
e_{j,j}d_{j,i}e_{i,i}=e_{j,j}d_{i,j}e_{i,i},\,\, i.e.\,\, d_{j,i}^{j,i}=d_{i,j}^{j,i}.
$$
From
$$
e_{i,i}[d_{i,j}\lambda e_{i,j}-\lambda e_{i,j}d_{i,j}+d_{j,i}\mu e_{j,i}-\mu e_{j,i}d_{j,i}]e_{j,j}=e_{i,i}[d(\lambda e_{i,j}+\mu e_{j,i})-(\lambda e_{i,j}+\mu e_{j,i})d]e_{j,j}
$$
we have that
$$
e_{i,i}d_{i,j}\lambda e_{i,j}-\lambda e_{i,j}d_{i,j}e_{j,j}=e_{i,i}d\lambda e_{i,j}-\lambda e_{i,j}de_{j,j},
$$
i.e.
$$
d_{i,j}^{i,i}-d_{i,j}^{j,j}=d^{i,i}-d^{j,j}.
$$
Similarly we have
$$
e_{j,j}d_{j,i}\mu e_{j,i}-\mu e_{j,i}d_{j,i}e_{i,i}=e_{j,j}d\mu e_{j,i}-\mu e_{j,i}de_{i,i},
$$
i.e.
$$
d_{j,i}^{j,j}-d_{j,i}^{i,i}=d^{j,j}-d^{i,i}.
$$
Therefore
$$
d_{i,j}^{i,i}-d_{i,j}^{j,j}=d_{j,i}^{i,i}-d_{j,i}^{j,j}.
$$

Now, let $e=e_{i,i}+e_{j,j}$. Then
$$
(1-e)d_{i,j}e_{j,j}e_{i,j}=e_{i,j}(1-e)d_{i,j}e_{j,j}=0, e_{j,j}d_{j,i}(1-e)e_{j,i}=e_{j,i}e_{j,j}d_{j,i}(1-e)=0,
$$
$$
(1-e)d_{j,i}e_{i,i}e_{j,i}=e_{j,i}(1-e)d_{j,i}e_{i,i}=0, e_{i,i}d_{i,j}(1-e)e_{i,j}=e_{i,j}e_{i,i}d_{i,j}(1-e)=0
$$
So we may assume that
$$
(1-e)d_{i,j}e_{j,j}=(1-e)d_{j,i}e_{j,j}, e_{j,j}d_{j,i}(1-e)=e_{j,j}d_{i,j}(1-e),
$$
$$
(1-e)d_{j,i}e_{i,i}=(1-e)d_{i,j}e_{i,i}, e_{i,i}d_{i,j}(1-e)=e_{i,i}d_{j,i}(1-e).
$$

So
$$
\nabla(\lambda e_{i,j})=e_{i,i}d_{i,j}\lambda e_{i,j}+e_{j,j}d_{i,j}\lambda e_{i,j}+(1-e)d_{i,j}\lambda e_{i,j}-\lambda e_{i,j}d_{i,j}e_{j,j}-
$$
$$
\lambda e_{i,j}d_{i,j}e_{i,i}-\lambda e_{i,j}d_{i,j}(1-e)=
$$
$$
e_{i,i}d_{i,j}\lambda e_{i,j}-\lambda e_{i,j}d_{i,j}e_{j,j}+e_{j,j}d_{i,j}e_{i,i}\lambda e_{i,j}-\lambda e_{i,j}e_{j,j}d_{i,j}e_{i,i}+
$$
$$
(1-e)d_{i,j}\lambda e_{i,j}-\lambda e_{i,j}d_{i,j}(1-e)=
$$
$$
e_{i,i}d_{j,i}\lambda e_{i,j}-\lambda e_{i,j}d_{j,i}e_{j,j}+e_{j,j}d_{j,i}e_{i,i}\lambda e_{i,j}-\lambda e_{i,j}e_{j,j}d_{j,i}e_{i,i}+
$$
$$
(1-e)d_{j,i}\lambda e_{i,j}-\lambda e_{i,j}d_{j,i}(1-e)=
$$
$$
d_{j,i}\lambda e_{i,j}-\lambda e_{i,j}d_{j,i}=D_{d_{j,i}}(\lambda e_{i,j}),
$$
since
$$
d_{i,j}^{i,j}=d_{j,i}^{i,j}, d_{i,j}^{j,i}=d_{j,i}^{j,i}.
$$
The proof is complete.
\end{proof}

\begin{lemma} \label{4}
For arbitrary $\lambda$, $\mu$, $\nu\in \mathcal{F}$ and each pair $i$, $j$ of distinct indices there exists an element $a\in M_n(\mathcal{F})$ such that
$$
\nabla(\lambda e_{i,i})=D_a(\lambda e_{i,i}), \nabla(\mu e_{i,j})=D_a(\mu e_{i,j}), \nabla(\nu e_{j,i})=D_a(\nu e_{j,i}).
$$
\end{lemma}

\begin{proof}
By lemmas \ref{1}, \ref{2} and \ref{3} there exist $a$, $b$, $c\in M_n(\mathcal{F})$ such that
$$
\nabla(\lambda e_{i,i})=D_a(\lambda e_{i,i}), \nabla(\mu e_{i,j})=D_a(\mu e_{i,j}),
$$
$$
\nabla(\lambda e_{i,i})=D_b(\lambda e_{i,i}), \nabla(\nu e_{j,i})=D_b(\nu e_{j,i})
$$
and
$$
\nabla(\mu e_{i,j})=D_c(\mu e_{i,j}), \nabla(\nu e_{j,i})=D_c(\nu e_{j,i}).
$$
We have
$$
D_a(\lambda e_{i,i})=D_b(\lambda e_{i,i}), D_a(\mu e_{i,j})=D_c(\mu e_{i,j})
$$
and
$$
D_b(\nu e_{j,i})=D_c(\nu e_{j,i}).
$$
Further from
$$
(D_a(\lambda e_{i,i}))e_{j,j}=(D_b(\lambda e_{i,i}))e_{j,j}, e_{j,j}(D_a(\lambda e_{i,i}))=e_{j,j}(D_b(\lambda e_{i,i}))
$$
we obtain that
$$
e_{i,i}ae_{j,j}=e_{i,i}be_{j,j}, e_{j,j}ae_{i,i}=e_{j,j}be_{i,i}.
$$
Similarly, from
$$
e_{i,i}(D_a(\mu e_{i,j}))e_{j,j}=e_{i,i}(D_c(\mu e_{i,j}))e_{j,j}, e_{j,j}(D_b(\nu e_{j,i}))e_{i,i}=e_{j,j}(D_c(\nu e_{j,i}))e_{i,i}
$$
it follows that
$$
e_{i,i}ae_{i,j}-e_{i,j}ae_{j,j}=e_{i,i}ce_{i,j}-e_{i,j}ce_{j,j}, e_{j,j}be_{j,i}-e_{j,i}be_{i,i}=e_{j,j}ce_{j,i}-e_{j,i}ce_{i,i}
$$
i.e.
$$
a^{i,i}-a^{j,j}=c^{i,i}-c^{j,j}, b^{j,j}-b^{i,i}=c^{j,j}-c^{i,i}.
$$
Hence
$$
a^{i,i}-a^{j,j}=b^{j,j}-b^{i,i}.
$$
Also
$$
(D_a(\lambda e_{i,i}))(1-e)=(D_b(\lambda e_{i,i}))(1-e)
$$
gives us $e_{i,i}a(1-e)=e_{i,i}b(1-e)$, and by the equality
$$
(1-e)D_b(\nu e_{j,i})=(1-e)D_c(\nu e_{j,i})
$$
we have $(1-e)be_{j,i}=(1-e)ce_{j,i}$.

At the same time $(1-e)ae_{j,j}e_{i,j}=e_{i,j}(1-e)ae_{j,j}=0$, $(1-e)ce_{j,j}e_{i,j}=e_{i,j}(1-e)ce_{j,j}=0$.

Hence we may assume $(1-e)be_{j,j}=(1-e)ce_{j,j}=(1-e)ae_{j,j}$.

Therefore,
$$
D_b(\nu e_{j,i})=b\nu e_{j,i}-\nu e_{j,i}b=e_{i,i}b\nu e_{j,i}+e_{j,j}b\nu e_{j,i}+(1-e)b\nu e_{j,i}-\nu e_{j,i}be_{i,i}-
$$
$$
\nu e_{j,i}be_{j,j}-\nu e_{j,i}b(1-e)=
$$
$$
e_{i,i}a\nu e_{j,i}+(b^{j,j}-b^{i,i})\nu e_{j,i}-\nu e_{j,i}ae_{j,j}+(1-e)a\nu e_{j,i}-\nu e_{j,i}a(1-e)=
$$
$$
e_{i,i}a\nu e_{j,i}+(a^{j,j}-a^{i,i})\nu e_{j,i}-\nu e_{j,i}ae_{j,j}+(1-e)a\nu e_{j,i}-\nu e_{j,i}a(1-e)=
$$
$$
a\nu e_{j,i}-\nu e_{j,i}a=D_a(\nu e_{j,i}).
$$
The proof is complete.
\end{proof}

\begin{lemma} \label{6}
For arbitrary $\lambda$, $\mu$, $\nu$ $\rho\in \mathcal{F}$ and each pair $i$, $j$ of distinct indices there exists an element $a\in M_n(\mathcal{F})$ such that
$$
\nabla(\lambda e_{i,i})=D_a(\lambda e_{i,i}), \nabla(\mu e_{i,j})=D_a(\mu e_{i,j}),
$$
$$
\nabla(\nu e_{j,i})=D_a(\nu e_{j,i}), \nabla(\rho e_{j,j})=D_a(\rho e_{j,j}).
$$
\end{lemma}

\begin{proof}
By lemma \ref{4} there exist $a$, $b\in M_n(\mathcal{F})$ such that
$$
\nabla(\lambda e_{i,i})=D_a(\lambda e_{i,i}), \nabla(\mu e_{i,j})=D_a(\mu e_{i,j}), \nabla(\nu e_{j,i})=D_a(\nu e_{j,i})
$$
$$
\nabla(\mu e_{i,j})=D_b(\mu e_{i,j}), \nabla(\nu e_{j,i})=D_b(\nu e_{j,i}), \nabla(\rho e_{j,j})=D_b(\rho e_{j,j}).
$$
We have
$$
D_a(\mu e_{i,j})=D_b(\mu e_{i,j}), D_a(\nu e_{j,i})=D_b(\nu e_{j,i}).
$$
From
$$
e_{j,j}(D_a(\mu e_{i,j}))=e_{j,j}(D_b(\mu e_{i,j})), e_{i,i}(D_a(\nu e_{j,i}))=e_{i,i}(D_b(\nu e_{j,i}))
$$
it follows that
$$
e_{j,j}ae_{i,i}=e_{j,j}be_{i,i}, e_{i,i}ae_{j,j}=e_{i,i}be_{j,j}
$$
respectively.
Also, by the equalities
$$
(D_a(\mu e_{i,j}))(1-e)=(D_b(\mu e_{i,j}))(1-e), (1-e)(D_a(\nu e_{j,i}))=(1-e)(D_b(\nu e_{j,i}))
$$
we have $e_{j,j}a(1-e)=e_{j,j}a(1-e)$ and $(1-e)ae_{j,j}=(1-e)be_{j,j}$ respectively,
where $e=e_{i,i}+e_{j,j}$.

Therefore
$$
D_b(\rho e_{j,j})=b\rho e_{j,j}-\rho e_{j,j}b=
$$
$$
e_{i,i}b\rho e_{j,j}+e_{j,j}b\rho e_{j,j}+(1-e)b\rho e_{j,j}-\rho e_{j,j}be_{i,i}-\rho e_{j,j}be_{j,j}-\rho e_{j,j}b(1-e)=
$$
$$
e_{i,i}b\rho e_{j,j}-\rho e_{j,j}be_{i,i}+(1-e)b\rho e_{j,j}-\rho e_{j,j}b(1-e)=
$$
$$
e_{i,i}a\rho e_{j,j}-\rho e_{j,j}ae_{i,i}+(1-e)a\rho e_{j,j}-\rho e_{j,j}a(1-e)=
$$
$$
e_{i,i}a\rho e_{j,j}-\rho e_{j,j}ae_{i,i}+e_{j,j}a\rho e_{j,j}-\rho e_{j,j}ae_{j,j}+(1-e)a\rho e_{j,j}-\rho e_{j,j}a(1-e)=
$$
$$
a\rho e_{j,j}-\rho e_{j,j}a=D_a(\rho e_{j,j}).
$$
This completes the proof.
\end{proof}

\begin{lemma} \label{7}
For arbitrary $\lambda$, $\mu\in \mathcal{F}$ and each index $i$ there exists an element $a\in M_n(\mathcal{F})$ such that
$$
\nabla(\lambda e_{i,i})=D_a(\lambda e_{i,i}), \nabla(\mu e_{i,i})=D_a(\mu e_{i,i}).
$$
\end{lemma}

\begin{proof}
By Lemma \ref{4} for arbitrary $\nu$, $\rho\in \mathcal{F}$ there exist $a$, $b\in M_n(\mathcal{F})$ such that
$$
\nabla(\lambda e_{i,i})=D_a(\lambda e_{i,i}), \nabla(\nu e_{i,j})=D_a(\nu e_{i,j}), \nabla(\rho e_{j,i})=D_a(\rho e_{j,i})
$$
$$
\nabla(\mu e_{i,i})=D_b(\mu e_{i,i}), \nabla(\nu e_{i,j})=D_b(\nu e_{i,j}), \nabla(\rho e_{j,i})=D_b(\rho e_{j,i}).
$$
We have
$$
D_a(\nu e_{i,j})=D_b(\nu e_{i,j}), D_a(\rho e_{j,i})=D_b(\rho e_{j,i}).
$$
Since
$$
e_{j,j}(D_a(\nu e_{i,j}))=e_{j,j}(D_b(\nu e_{i,j})), e_{i,i}(D_a(\rho e_{j,i}))=e_{i,i}(D_b(\rho e_{j,i}))
$$
we have that
$$
e_{j,j}ae_{i,i}=e_{j,j}be_{i,i}, e_{i,i}ae_{j,j}=e_{i,i}be_{j,j}
$$
respectively.
Let $e=e_{i,i}+e_{j,j}$. Then by the equalities
$$
(1-e)(D_a(\nu e_{i,j}))=(1-e)(D_b(\nu e_{i,j})), (D_a(\rho e_{j,i}))(1-e)=(D_b(\rho e_{j,i}))(1-e)
$$
we have $(1-e)ae_{i,i}=(1-e)be_{i,i}$ and $e_{i,i}a(1-e)=e_{i,i}b(1-e)$ respectively.

Therefore
$$
D_b(\mu e_{i,i})=b\mu e_{i,i}-\mu e_{i,i}b=
$$
$$
e_{i,i}b\mu e_{i,i}+e_{j,j}b\mu e_{i,i}+(1-e)b\mu e_{i,i}-\mu e_{i,i}be_{i,i}-\mu e_{i,i}be_{j,j}-\mu e_{i,i}b(1-e)=
$$
$$
e_{j,j}b\mu e_{i,i}-\mu e_{i,i}be_{j,j}+(1-e)b\mu e_{i,i}-\mu e_{i,i}b(1-e)=
$$
$$
e_{j,j}a\mu e_{i,i}-\mu e_{i,i}ae_{j,j}+(1-e)a\mu e_{i,i}-\mu e_{i,i}a(1-e)=
$$
$$
e_{j,j}a\mu e_{i,i}-\mu e_{i,i}ae_{j,j}+e_{i,i}a\mu e_{i,i}-\mu e_{i,i}ae_{i,i}+(1-e)a\mu e_{i,i}-\mu e_{i,i}a(1-e)=
$$
$$
a\mu e_{i,i}-\mu e_{i,i}a=D_a(\mu e_{i,i}).
$$
This completes the proof.
\end{proof}

Similarly we can prove the following lemma using the above lemmas \ref{4} and \ref{6}.

\begin{lemma} \label{8}
For arbitrary $\lambda$, $\mu\in \mathcal{F}$ and each pair $i$, $j$ of distinct indices there exist elements $a\in M_n(\mathcal{F})$ such that
$$
\nabla(\lambda e_{i,j})=D_a(\lambda e_{i,j}), \nabla(\mu e_{i,j})=D_a(\mu e_{i,j}).
$$
\end{lemma}

\begin{lemma} \label{9}
For arbitrary $\lambda$, $\mu\in \mathcal{F}$ and each pair $i$, $j$ of distinct indices one has
$$
\nabla(\lambda e_{i,i})=\lambda \nabla(e_{i,i}), \nabla(\mu e_{i,j})=\mu \nabla(e_{i,j}).
$$
\end{lemma}

\begin{proof}
By lemma \ref{7} there exists an element $a\in M_n(\mathcal{F})$ such that
$$
\nabla(\lambda e_{i,i})=D_a(\lambda e_{i,i}), \nabla(e_{i,i})=D_a(e_{i,i}).
$$
Hence
$$
\nabla(\lambda e_{i,i})=D_a(\lambda e_{i,i})=\lambda D_a(e_{i,i})=\lambda\nabla(e_{i,i}).
$$
The second equality is  proved in a similarly way.
\end{proof}

\begin{theorem} \label{141}
Let $\mathcal{F}$ be an arbitrary field, and let $M_2(\mathcal{F})$ be the algebra
of $2\times 2$ matrices over $\mathcal{F}$. Then any additive local inner derivation on the matrix algebra
$M_2(\mathcal{F})$ is an inner derivation.
\end{theorem}

\begin{proof}
Let $\nabla :M_2(\mathcal{F})\to M_2(\mathcal{F})$ be an additive local inner derivation.
By lemma \ref{6} there exists $a\in M_2(\mathcal{F})$ such that
$$
\nabla (e_{i,j})=D_a(e_{i,j})
$$
for any indices $i$, $j$ from $\{1,2\}.$

Let $x$ be an arbitrary element in $M_2(\mathcal{F})$.
Then $x=\sum_{k,l=1}^2 x^{k,l}e_{k,l}$ and by lemma \ref{9}
$$
\nabla(x)=\sum_{k,l=1}^2 \nabla(x^{k,l}e_{k,l})=\sum_{k,l=1}^2 x^{k,l}\nabla(e_{k,l})=
$$
$$
\sum_{k,l=1}^2 x^{k,l}D_a(e_{k,l})=D_a(x).
$$
Hence $\nabla$ is an inner derivation. The proof is complete.
\end{proof}

\begin{lemma} \label{10}
Let $\nabla :M_n(\mathcal{F})\to M_n(\mathcal{F})$ be an additive local inner derivation.
Then for any indices $i$, $j$, $k$, $l$, at least three of which are pairwise distinct,
there exists $a\in M_n(\mathcal{F})$ such that
$$
\nabla (e_{i,j})=D_a(e_{i,j}), \nabla (e_{k,l})=D_a(e_{k,l}).
$$
\end{lemma}

\begin{proof}
Let $i$, $j$, $k$ be pairwise distinct indices. We have
$$
\nabla(\lambda e_{i,j})+\nabla(\mu e_{j,k})=\nabla(\lambda e_{i,j}+\mu e_{j,k})
$$
Let $d_{i,j}$, $d_{j,k}$, $d\in M_n(\mathcal{F})$ be such elements that
$$
\nabla(\lambda e_{i,j})=D_{d_{i,j}}(\lambda e_{i,j}),
\nabla(\mu e_{j,k})=D_{d_{j,k}}(\mu e_{j,k}),
$$
$$
\nabla(\lambda e_{i,j}+\mu e_{j,k})=D_d(\lambda e_{i,j}+\mu e_{j,k}).
$$
Then
$$
D_{d_{i,j}}(\lambda e_{i,j})+D_{d_{j,k}}(\mu e_{j,k})=D_d(\lambda e_{i,j}+\mu e_{j,k}).
$$
Note that
$$
e_{j,j}d_{j,k}e_{i,i}\mu e_{j,k}=\mu e_{j,k}e_{j,j}d_{j,k}e_{i,i}=0.
$$
So we can take
$$
e_{j,j}d_{j,k}e_{i,i}=e_{j,j}d_{i,j}e_{i,i},\,\, i.e.\,\, d_{j,k}^{j,i}=d_{i,j}^{j,i}.
$$
 From
$$
e_{i,i}[D_{d_{i,j}}(\lambda e_{i,j})+D_{d_{j,k}}(\mu e_{j,k})]e_{j,j}=e_{i,i}[D_d(\lambda e_{i,j}+\mu e_{j,k})]e_{j,j}
$$
it follows that
$$
e_{i,i}d_{i,j}\lambda e_{i,j}-\lambda e_{i,j}d_{i,j}e_{j,j}=e_{i,i}d\lambda e_{i,j}-\lambda e_{i,j}de_{j,j},
$$
i.e.
$$
d_{i,j}^{i,i}-d_{i,j}^{j,j}=d^{i,i}-d^{j,j}.
$$
Similarly we have
$$
e_{j,j}d_{j,k}\mu e_{j,k}-\mu e_{j,k}d_{j,k}e_{i,i}=e_{j,j}d\mu e_{j,k}-\mu e_{j,k}de_{i,i},
$$
i.e.
$$
d_{j,k}^{j,j}-d_{j,k}^{i,i}=d^{j,j}-d^{i,i}.
$$
Therefore
$$
d_{i,j}^{i,i}-d_{i,j}^{j,j}=d_{j,k}^{i,i}-d_{j,k}^{j,j}.
$$

Now, let $e=e_{i,i}+e_{j,j}$. Then
$$
e_{j,j}d_{j,k}(1-e)e_{j,k}=e_{j,k}e_{j,j}d_{j,k}(1-e)=0, (1-e)d_{j,k}e_{i,i}e_{j,k}=e_{j,k}(1-e)d_{j,k}e_{i,i}=0.
$$
So we may assume that
$$
e_{j,j}d_{j,k}(1-e)=e_{j,j}d_{i,j}(1-e), (1-e)d_{j,k}e_{i,i}=(1-e)d_{i,j}e_{i,i}.
$$

So
$$
\nabla(\lambda e_{i,j})=e_{i,i}d_{i,j}\lambda e_{i,j}+e_{j,j}d_{i,j}\lambda e_{i,j}+(1-e)d_{i,j}\lambda e_{i,j}-\lambda e_{i,j}d_{i,j}e_{j,j}-
$$
$$
\lambda e_{i,j}d_{i,j}e_{i,i}-\lambda e_{i,j}d_{i,j}(1-e)=
$$
$$
e_{i,i}d_{i,j}\lambda e_{i,j}-\lambda e_{i,j}d_{i,j}e_{j,j}+e_{j,j}d_{i,j}e_{i,i}\lambda e_{i,j}-\lambda e_{i,j}e_{j,j}d_{i,j}e_{i,i}+
$$
$$
(1-e)d_{i,j}\lambda e_{i,j}-\lambda e_{i,j}d_{i,j}(1-e)=
$$
$$
e_{i,i}d_{j,k}\lambda e_{i,j}-\lambda e_{i,j}d_{j,k}e_{j,j}+e_{j,j}d_{j,k}e_{i,i}\lambda e_{i,j}-\lambda e_{i,j}e_{j,j}d_{j,k}e_{i,i}+
$$
$$
(1-e)d_{j,k}\lambda e_{i,j}-\lambda e_{i,j}d_{j,k}(1-e)=
$$
$$
d_{j,k}\lambda e_{i,j}-\lambda e_{i,j}d_{j,k}=D_{d_{j,k}}(\lambda e_{i,j}),
$$
since
$$
d_{i,j}^{j,i}=d_{j,k}^{j,i}.
$$

Let $i$, $j$, $k$, $l$ be pairwise distinct indices and let
$d_{i,j}$, $d_{k,l}$ be elements in $M_n(\mathcal{F})$ such that
$$
\nabla (e_{i,j})=D_{d_{i,j}}(e_{i,j}), \nabla (e_{k,l})=D_{d_{k,l}}(e_{k,l}).
$$
Put $e=e_{i,i}+e_{j,j}$, $f=e_{k,k}+e_{l,l}$.

Then we have
$$
(1-f)d_{k,l}(1-f)e_{k,l}=e_{k,l}(1-f)d_{k,l}(1-f)=0.
$$
So we may take $(1-f)d_{k,l}(1-f)=(1-f)d_{i,j}(1-f)$. Hence
$$
(e_{i,i}+e_{j,j})d_{i,j}(e_{i,i}+e_{j,j})=(e_{i,i}+e_{j,j})d_{k,l}(e_{i,i}+e_{j,j})
$$
and
$$
e_{i,i}d_{i,j}e_{i,i}=e_{i,i}d_{k,l}e_{i,i}, e_{i,i}d_{i,j}e_{j,j}=e_{i,i}d_{k,l}e_{j,j},
$$
$$
e_{j,j}d_{i,j}e_{i,i}=e_{j,j}d_{k,l}e_{i,i}, e_{j,j}d_{i,j}e_{j,j}=e_{j,j}d_{k,l}e_{j,j}.
$$
Also we have
$$
(1-e-e_{l,l})d_{k,l}e_{i,i}e_{k,l}=e_{k,l}(1-e-e_{l,l})d_{k,l}e_{i,i}=0,
$$
$$
e_{j,j}d_{k,l}(1-e-e_{k,k})e_{k,l}=e_{k,l}e_{j,j}d_{k,l}(1-e-e_{k,k})=0.
$$
So we may take
$$
(1-e-e_{l,l})d_{k,l}e_{i,i}=(1-e-e_{l,l})d_{i,j}e_{i,i},
$$
$$
e_{j,j}d_{k,l}(1-e-e_{k,k})=e_{j,j}d_{i,j}(1-e-e_{k,k}).
$$
Now, let $d$ be an element in $M_n(\mathcal{F})$ such that
$$
\nabla (e_{i,j}+e_{k,l})=D_d(e_{i,j}+e_{k,l}).
$$
Then $D_d(e_{i,j}+e_{k,l})=D_{d_{i,j}}(e_{i,j})+D_{d_{k,l}}(e_{k,l})$ and by the equalities
$$
e_{i,i}(D_d(e_{i,j}+e_{k,l}))e_{k,k}=e_{i,i}(D_{d_{i,j}}(e_{i,j})+D_{d_{k,l}}(e_{k,l}))e_{k,k},
$$
$$
e_{j,j}(D_d(e_{i,j}+e_{k,l}))e_{l,l}=e_{j,j}(D_{d_{i,j}}(e_{i,j})+D_{d_{k,l}}(e_{k,l}))e_{l,l}
$$
we have
$$
e_{j,j}de_{k,k}=e_{j,j}d_{i,j}e_{k,k}, e_{j,j}de_{k,k}=e_{j,j}d_{k,l}e_{k,k}.
$$
Hence $e_{j,j}d_{i,j}e_{k,k}=e_{j,j}d_{k,l}e_{k,k}$.

Similarly by the equalities
$$
e_{k,k}(D_d(e_{i,j}+e_{k,l}))e_{i,i}=e_{k,k}(D_{d_{i,j}}(e_{i,j})+D_{d_{k,l}}(e_{k,l}))e_{i,i},
$$
$$
e_{l,l}(D_d(e_{i,j}+e_{k,l}))e_{j,j}=e_{l,l}(D_{d_{i,j}}(e_{i,j})+D_{d_{k,l}}(e_{k,l}))e_{j,j}
$$
we have
$$
e_{l,l}d_{i,j}e_{i,i}=e_{l,l}d_{k,l}e_{i,i}.
$$
Therefore
$$
\nabla(e_{i,j})=e_{i,i}d_{i,j}e_{i,j}+e_{j,j}d_{i,j}e_{i,j}+(1-e)d_{i,j}e_{i,j}-e_{i,j}d_{i,j}e_{j,j}-
$$
$$
e_{i,j}d_{i,j}e_{i,i}-e_{i,j}d_{i,j}(1-e)=
$$
$$
e_{i,i}d_{i,j}e_{i,j}-e_{i,j}d_{i,j}e_{j,j}+e_{j,j}d_{i,j}e_{i,i}e_{i,j}-e_{i,j}e_{j,j}d_{i,j}e_{i,i}+
$$
$$
(1-e)d_{i,j}e_{i,j}-e_{i,j}d_{i,j}(1-e)=
$$
$$
e_{i,i}d_{k,l}e_{i,j}-e_{i,j}d_{k,l}e_{j,j}+e_{j,j}d_{k,l}e_{i,i}e_{i,j}-e_{i,j}e_{j,j}d_{k,l}e_{i,i}+
$$
$$
(1-e)d_{k,l}e_{i,j}-e_{i,j}d_{k,l}(1-e)=
$$
$$
d_{k,l}e_{i,j}-e_{i,j}d_{k,l}=D_{d_{k,l}}(e_{i,j}).
$$

Similarly we can prove that
for any pairwise distinct indices $i$, $j$, $k$ there exists
$a$ in $M_n(\mathcal{F})$ such that
$$
\nabla (e_{i,j})=D_a(e_{i,j}), \nabla (e_{k,k})=D_a(e_{k,k}).
$$
This completes the proof.
\end{proof}

\begin{lemma} \label{11}
There exists $a\in M_n(\mathcal{F})$ such that
$$
\nabla (e_{i,i+1})=D_a(e_{i,i+1}), i=1,2,\dots n-1.
$$
\end{lemma}

\begin{proof}
By lemma \ref{10} there exists $a\in M_n(\mathcal{F})$ such that $\nabla (e_{1,2})=D_a(e_{1,2})$, $\nabla (e_{2,3})=D_a(e_{2,3})$.
Suppose that for $k$ there exists $a\in M_n(\mathcal{F})$ such that
$$
\nabla (e_{i,i+1})=D_a(e_{i,i+1}), i=1,2,\dots k-1.
$$
We prove that for $k+1$ there exists $b\in M_n(\mathcal{F})$ such that
$$
\nabla (e_{i,i+1})=D_b(e_{i,i+1}), i=1,2,\dots k.
$$
By lemma \ref{10} there exists $c(i,i+1)\in M_n(\mathcal{F})$ such that $\nabla (e_{i,i+1})=D_{c(i,i+1)}(e_{i,i+1})$, $\nabla
(e_{k,k+1})=D_{c(i,i+1)}(e_{k,k+1})$,
where $i=1,2,\dots k-1$.
Since
$$
e_{k+1,k+1}ae_{k+1,k+1}e_{i,j}=e_{i,j}e_{k+1,k+1}ae_{k+1,k+1}=0, i,j=1,2,\dots k
$$
we may put
$$
b^{i,j}=a^{i,j},\,\, if\,\,\, i\leq k \,\,\,or\,\,\, j\leq k,
$$
$$
b^{k+1,k+1}=c(k,k+1)^{k+1,k+1}-c(k,k+1)^{k,k}+b^{k,k},
$$
and
$$
b^{i,j}=c(k,k+1)^{i,j}, (i,j)\neq (k+1,k+1), i\geq k+1, j\geq k+1.  \eqno{(5.1)}
$$
In this case we have
$$
e_{k,k}\nabla(e_{k,k+1})e_{k+1,k+1}=e_{k,k}D_{c(k,k+1)}(e_{k,k+1})e_{k+1,k+1}=
$$
$$
e_{k,k}c(k,k+1)e_{k,k+1}e_{k+1,k+1}-e_{k,k}e_{k,k+1}c(k,k+1)e_{k+1,k+1}=
$$
$$
e_{k,k}be_{k,k+1}e_{k+1,k+1}-e_{k,k}e_{k,k+1}be_{k+1,k+1}=e_{k,k}D_b(e_{k,k+1})e_{k+1,k+1}
$$
by equalities (5.1). If $i\neq k$ and $j\neq k+1$ then
$$
e_{i,i}\nabla(e_{k,k+1})e_{j,j}=0=e_{i,i}D_b(e_{k,k+1})e_{j,j} \eqno{(5.2)}
$$
If $j>k+1$ then
$$
e_{k,k}\nabla(e_{k,k+1})e_{j,j}=e_{k,k}D_{c(k,k+1)}(e_{k,k+1})e_{j,j}=
$$
$$
e_{k,k}c(k,k+1)e_{k,k+1}e_{j,j}-e_{k,k}e_{k,k+1}c(k,k+1)e_{j,j}=
$$
$$
e_{k,k}be_{k,k+1}e_{j,j}-e_{k,k}e_{k,k+1}be_{j,j}=e_{k,k}D_b(e_{k,k+1})e_{j,j}
$$
by (5.1). From
$$
c(j,j+1)e_{j,j+1}-e_{j,j+1}c(j,j+1)=ae_{j,j+1}-e_{j,j+1}a=be_{j,j+1}-e_{j,j+1}b
$$
it follows that
$$
e_{k+1,k+1}c(j,j+1)e_{j,j}=e_{k+1,k+1}be_{j,j},          \eqno{(5.3)}
$$
where $j=1,2,\dots k$. Therefore, if $j<k+1$ then
$$
e_{k,k}\nabla(e_{k,k+1})e_{j,j}=e_{k,k}D_{c(j,j+1)}(e_{k,k+1})e_{j,j}=
$$
$$
e_{k,k}c(j,j+1)e_{k,k+1}e_{j,j}-e_{k,k}e_{k,k+1}c(j,j+1)e_{j,j}=
$$
$$
e_{k,k}be_{k,k+1}e_{j,j}-e_{k,k}e_{k,k+1}be_{j,j}=e_{k,k}D_b(e_{k,k+1})e_{j,j}
$$
by (5.3). Similarly
$$
e_{j,j}\nabla(e_{k,k+1})e_{k+1,k+1}=e_{j,j}D_b(e_{k,k+1})e_{k+1,k+1}, j=1,2,\dots n.
$$
So
$$
\nabla(e_{k,k+1})=\sum_{i,j=1}^n e_{i,i}\nabla(e_{k,k+1})e_{j,j}=\sum_{i,j=1}^n e_{i,i}D_b(e_{k,k+1})e_{j,j}=D_b(e_{k,k+1}).
$$
Hence by the induction we obtain that there exists $a\in M_n(\mathcal{F})$ such that
$$
\nabla (e_{i,i+1})=D_a(e_{i,i+1}), i=1,2,\dots n-1.
$$
\end{proof}

\begin{lemma} \label{12}
For any indices $i$, $j$ there exists $a\in M_n(\mathcal{F})$ such that
$$
\nabla(\sum_{k=1}^{n-1} e_{k,k+1})=D_a(\sum_{k=1}^{n-1} e_{k,k+1}),
\nabla(e_{i,j})=D_a(e_{i,j}).
$$
\end{lemma}

\begin{proof} By lemma \ref{11} there exists $a\in M_n(\mathcal{F})$ such that
$$
\nabla (e_{i,i+1})=D_a(e_{i,i+1}), i=1,2,\dots n-1.
$$

Fix $i$ and $j$ from $\{1,2,\dots n\}$ such that $i\leq j$.
We have $\nabla(\sum_{k=i}^{j-1} e_{k,k+1}+e_{j,i})=\nabla(\sum_{k=i}^{j-1} e_{k,k+1})+\nabla(e_{j,i})$.
There exist $b$, $d\in M_n(\mathcal{F})$ such that $\nabla(\sum_{k=i}^{j-1} e_{k,k+1}+e_{j,i})=D_d(\sum_{k=i}^{j-1} e_{k,k+1}+e_{j,i})$,
$\nabla(e_{j,i})=D_b{e_{j,i}}$. Therefore
$$
D_d(\sum_{k=i}^{j-1}e_{k,k+1}+e_{j,i})=D_a(\sum_{k=i}^{j-1}e_{k,k+1})+D_b(e_{j,i}).
$$
Hence
$$
d^{i,i}-d^{i+1,i+1}=a^{i,i}-a^{i+1,i+1}, k=i,i+1,i+2,\dots,j,  d^{j,j}-d^{i,i}=b^{j,j}-b^{i,i},
$$
and
$$
a^{j,j}-a^{i,i}=b^{j,j}-b^{i,i}.
$$

Now, if $k\neq i$ and $k\neq j$ then
$$
e_{k,k}be_{k,k}e_{j,i}=e_{j,i}e_{k,k}be_{k,k}=0.
$$
So we may assume that $e_{k,k}be_{k,k}=e_{k,k}ae_{k,k}$ for every $k$ such that $k\neq i$ and $k\neq j$.
By lemma \ref{10} there exists $c(i,i+1)\in M_n(\mathcal{F})$ such that $\nabla (e_{i,i+1})=D_{c(i,i+1)}(e_{i,i+1})$, $\nabla (e_{j,i})=D_{c(i,i+1)}(e_{j,i})$.
We have $D_{c(i,i+1)}(e_{i,i+1})=D_a(e_{i,i+1})$, $D_{c(i,i+1)}(e_{j,i})=D_b(e_{j,i})$
and
$$
e_{j,j}c(i,i+1)e_{i,i+1}=e_{j,j}ae_{i,i+1}, e_{j,j}be_{i,i}e_{j,i}=e_{j,i}e_{j,j}be_{i,i}=0.
$$
So we may assume that $e_{j,j}be_{i,i}=e_{j,j}c(i,i+1)e_{i,i}$. Hence $e_{j,j}be_{i,i}=e_{j,j}ae_{i,i}$. Also, from
$$
D_{c(i-1,i)}(e_{i-1,i})=D_a(e_{i-1,i}), D_{c(i-1,i)}(e_{j,i})=D_b(e_{j,i})
$$
it follows that
$$
e_{i-1,i}c(i-1,i)e_{j,j}=e_{i-1,i}ae_{j,j}, e_{i,i}c(i-1,i)e_{j,i}=e_{i,i}be_{j,i},
$$
and
$$
e_{i,i}be_{j,j}=e_{i,i}ae_{j,j}.
$$

The  remaining case is $\{k,l\}\neq \{i,j\}$ for $e_{k,l}$.

1) Suppose $k=i$, $l\neq j$, $l\neq i$. Take
$$
D_{c(i-1,i)}(e_{i-1,i})=D_a(e_{i-1,i}), D_{c(i-1,i)}(e_{j,i})=D_b(e_{j,i})
$$
and we have
$$
e_{i-1,i}c(i-1,i)e_{l,l}=e_{i-1,i}ae_{l,l}, e_{j,i}c(i-1,i)e_{l,l}=e_{j,i}be_{l,l},
$$
and so $e_{i,i}be_{l,l}=e_{i,i}ae_{l,l}$.

Now, we take $D_{c(i,i+1)}(e_{i,i+1})=D_a(e_{i,i+1})$, $D_{c(i,i+1)}(e_{j,i})=D_b(e_{j,i})$. Then
$$
e_{l,l}c(i,i+1)e_{i,i+1}=e_{l,l}ae_{i,i+1}, e_{l,l}be_{i,i}e_{j,i}=e_{j,i}e_{l,l}be_{i,i}=0.
$$
So we may assume $e_{l,l}be_{i,i}=e_{l,l}ae_{i,i}$.

2) Suppose $k\neq i$, $l=j$, $k\neq j$. Take
$$
D_{c(j-1,j)}(e_{j-1,j})=D_a(e_{j-1,j}), D_{c(j-1,j)}(e_{j,i})=D_b(e_{j,i})
$$
and we have
$$
e_{j-1,j}c(j-1,j)e_{k,k}=e_{j-1,j}ae_{k,k}, e_{j,j}be_{k,k}e_{j,i}=e_{j,i}e_{j,j}be_{k,k}=0.
$$
So we may assume that $e_{j,j}be_{k,k}=e_{j,j}ae_{k,k}$.

Take
$$
D_{c(j,j+1)}(e_{j,j+1})=D_a(e_{j,j+1}), D_{c(j,j+1)}(e_{j,i})=D_b(e_{j,i})
$$
and we have
$$
e_{k,k}c(j,j+1)e_{j,j+1}=e_{k,k}ae_{j,j+1}, e_{k,k}c(j,j+1)e_{j,i}=e_{k,k}be_{j,i},
$$
and so $e_{k,k}be_{j,j}=e_{k,k}ae_{j,j}$.

3) Now, suppose $k\neq i,j$, $l\neq i,j$. Then
$e_{k,k}be_{l,l}e_{j,i}=e_{j,i}e_{k,k}be_{l,l}=0$.
So we may assume that $e_{k,k}be_{l,l}=e_{k,k}ae_{l,l}$.
Thus, for all $\{k,l\}\neq \{i,j\}$ we have $e_{k,k}be_{l,l}=e_{k,k}ae_{l,l}$ and,
if $\{k,l\}=\{i,j\}$ then $b^{i,i}-b^{j,j}=a^{i,i}-a^{j,j}$ and $b^{i,j}=b^{i,j}$,
$b^{j,i}=b^{j,i}$. Hence
$$
\nabla (e_{j,i})=D_b(e_{j,i})=D_a(e_{j,i}).
$$

Now, by the definition of additive local inner derivation we have that  $\nabla (e_{i,j})=D_c(e_{i,j})$, $\nabla (e_{j,i})=D_c(e_{j,i})$ for some $c$ in
$M_n(\mathcal{F})$. Then $D_c(e_{j,i})=D_a(e_{j,i})$ and $c^{j,j}-c^{i,i}=a^{j,j}-a^{i,i}$,
$e_{i,i}ce_{j,j}=e_{i,i}ae_{j,j}$. Also we have
$$
e_{j,j}ce_{i,i}e_{j,i}=e_{j,i}e_{j,j}ce_{i,i}=0.
$$
So we may assume that $e_{j,j}ce_{i,i}=e_{j,j}ae_{i,i}$.

Now similar to the equality $\nabla (e_{j,i})=D_a(e_{j,i})$ we prove that  $\nabla (e_{i,j})=D_c(e_{i,j})=D_a(e_{i,j})$.
This completes the proof.
\end{proof}

\begin{lemma} \label{13}
There exists $a\in M_n(\mathcal{F})$ such that
$$
\nabla (e_{i,j})=D_a(e_{i,j})
$$
for any indices $i$, $j$.
\end{lemma}

\begin{proof}
By the previous lemmas we can repeat the proof of theorem 4 in \cite{AA3}
and get the statement of this lemma. The proof is complete.
\end{proof}

\begin{theorem} \label{14}
Let $\mathcal{F}$ be an arbitrary field, and let $M_n(\mathcal{F})$ be the algebra
of $n\times n$ matrices over $\mathcal{F}$, $n>1$. Then any additive local inner derivation on the algebra
$M_n(\mathcal{F})$ is an inner derivation.
\end{theorem}

\begin{proof}
Let $\nabla :M_n(\mathcal{F})\to M_n(\mathcal{F})$ be an additive local inner derivation. Then
by lemma \ref{13} there exists $a\in M_n(\mathcal{F})$ such that for any indices $i$, $j$
$$
\nabla (e_{i,j})=D_a(e_{i,j}).
$$
Let $x$ be an arbitrary element in $M_n(\mathcal{F})$.
Then $x=\sum_{k,l=1}^n x^{k,l}e_{k,l}$ and by lemma \ref{9} we have that
$$
\nabla (x)=\sum_{k,l=1}^n \nabla(x^{k,l}e_{k,l})=\sum_{k,l=1}^n x^{k,l}\nabla(e_{k,l})=
$$
$$
\sum_{k,l=1}^n x^{k,l}D_a(e_{k,l})=D_a(x).
$$
Thus $\nabla$ is an inner derivation. The proof is complete.
\end{proof}

Let $\mathcal{R}$ be a ring. An additive map $D : \mathcal{R}\to \mathcal{R}$ is called
a derivation, if $D(xy)=D(x)y+xD(y)$ for any two elements $x$, $y\in \mathcal{R}$.

An additive map $\nabla : \mathcal{R}\to \mathcal{R}$ is called local derivation, if for
any element $x\in \mathcal{R}$ there exists a derivation
$D:\mathcal{R}\to \mathcal{R}$ such that $\nabla (x)=D(x)$.

Now let $\mathcal{R}$ be a non commutative (but associative) ring. A derivation $D$ on $\mathcal{R}$
is called an inner derivation, if there exists an element $a\in \mathcal{R}$ such that
$$
D(x)=ax-xa, x\in \mathcal{R}.
$$
This derivation $D$ we denote by $D_a$, i.e. $D_a(x)=ax-xa$.
An additive map $\nabla : \mathcal{R}\to \mathcal{R}$ is called local inner derivation,
if for any element $x\in \mathcal{R}$ there exists an inner derivation $D_a$ such that
$\nabla(x)=D_a(x)$.

A finite ring is a ring that has a finite number of elements. A finite ring generated by the identity element
is a finite ring, every element of which is a sum of some quantity of the identity element of this ring.
By the proofs of the previous lemmas we have the following lemma.

\begin{lemma} \label{15}
Let $\Re$ be a finite ring generated by the identity element
or the ring of integers, $M_n(\Re)$ be the ring
of $n\times n$ matrices over $\Re$, $n>1$, and let $\nabla :M_n(\Re)\to M_n(\Re)$ be a local inner derivation. Then
there exists $a\in M_n(\lambda\Re)$ such that
$$
\nabla (e_{i,j})=D_a(e_{i,j})
$$
for any indices $i$, $j$.
\end{lemma}

By lemma \ref{15} and by the definition of a local inner derivation we have the following theorem.

\begin{theorem} \label{16}
Let $\Re$ be a finite ring generated by the identity element
or the ring of integers, $M_n(\Re)$ be the ring
of $n\times n$ matrices over $\Re$, $n>1$. Then for every local inner derivation $\nabla$ on
the ring $M_n(\Re)$ there exists a matrix $a\in M_n(\Re)$ such that
$$
\nabla(x)=D_a(x), x\in M_n(\Re),
$$
i.e. $\nabla$ is a derivation.
\end{theorem}

Let $\mathcal{B}$ be a subalgebra of an algebra $\mathcal{A}$. A derivation $D$
on $\mathcal{B}$ is said to be spatial, if $D$ is implemented by an element in $\mathcal{A}$, i.e.
$$
D(x)=ax-xa, x\in \mathcal{B},
$$
for some $a\in \mathcal{A}$. An additive local derivation $\nabla$ on $\mathcal{B}$ is called additive local spatial derivation with respect
to derivations implemented by an element in $\mathcal{A}$, if for every element $x\in \mathcal{B}$ there
exists an element $a\in \mathcal{A}$ such that $\nabla(x)=ax-xa$.

It should be noted that by the proofs of theorems 5.4, 5.5 and 6.1 in \cite{AA4} the
notions of local inner derivation and local spatial derivation in these theorems can be replaced by the
notions of additive local inner derivation and additive local spatial derivation respectively.

\section{Local derivations on Jordan algebras of symmetric matrices}

This section is devoted to derivations and local derivations of Jordan algebras.
In this section the  notations and terminology  follow the paper \cite{UH}
of H. Upmeier.

Given subsets $B$ and $C$ of a Lie ring $\Re$ with bracket $[\cdot, \cdot]$, let $[B,C]$ denote the subset of $\Re$
consisting of all finite sums of elements $[b,c]$, where $b\in B$ and $c\in C$.

Consider a Jordan ring $\mathcal{J}$ and let $m=\{xM: x\in \mathcal{A}\}$, where $xM$  denotes the multiplication
operator defined by $(xM)y:=x\cdot y$ for all  $y\in \mathcal{A}$.
Let $aut(\mathcal{J})$ denotes the Lie algebra of all derivations of $\mathcal{J}$.
The elements of the ideal $int(\mathcal{J}):=[m,m]$  in $aut(\mathcal{J})$ are called $inner$ derivations of the Jordan ring $\mathcal{J}$ (cf.\cite{UH}).

Let $\mathcal{R}$ be an associative unital ring, and suppose $2$ is invertible in $\mathcal{R}$
Then the set $\mathcal{R}$ with respect to the operations of
addition and Jordan multiplication
$$
a\cdot b=\frac{1}{2}(ab+ba), a, b\in \mathcal{R}
$$
is a Jordan ring. This Jordan ring we will denote by $(\mathcal{R}, \cdot)$.
Every inner derivation of $(\mathcal{R}, \cdot)$ is an inner derivation of $\mathcal{R}$, and,
conversely, every inner derivation of $\mathcal{R}$ is an inner derivation of $(\mathcal{R}, \cdot)$ \cite{AA3}.

Let $\nabla$ be a local inner derivation of the Jordan ring $(\mathcal{R}, \cdot)$. Then for every element $x\in \mathcal{R}$
there is an inner derivation $D$ of $(\mathcal{R}, \cdot)$ such that $\nabla(x)=D(x)$.
But $D$ is also an inner derivation of the associative ring $\mathcal{R}$. Hence, $\nabla$ is a local inner derivation
of the associative ring $\mathcal{R}$. Conversely, every local inner derivation of the associative ring $\mathcal{R}$
is a local inner derivation of the Jordan ring $(\mathcal{R}, \cdot)$.

Now, let $\mathcal{R}$ be an involutive unital ring, and suppose $2$ is invertible in $\mathcal{R}$. Let $\mathcal{R}_ {sa}$ be the set of all self-adjoint elements of
the ring $\mathcal{R}$. Then, it is known that $(\mathcal{R} _{sa}, \cdot)$ is a Jordan ring.
Also, every inner derivation of the Jordan ring $(\mathcal{R}_{sa}, \cdot)$ is extended to
an inner derivation of the $*$-ring $\mathcal{R}$ \cite{AA3}. Such extension of derivations on
a special Jordan algebra is considered in \cite{UH}. Concerning local inner derivation, till now it is not possible
to obtain such extension without additional conditions. This problem shows the importance of the main result
in the present section.

Throughout of this section  $\mathcal{F}$ is an arbitrary field with invertible $2$, and $M_n(\mathcal{F})$ is the associative algebra of $n\times n$ matrices
over $\mathcal{F}$. In this case the set
$$
H_n(\mathcal{F})=\{
\left[%
\begin{array}{cccc}
a^{1,1} & a^{1,2} & \cdots & a^{1,n}\\
a^{2,1} & a^{2,2} & \cdots & a^{2,n}\\
\vdots & \vdots & \ddots & \vdots\\
a^{n,1} & a^{n,2} & \cdots & a^{n,n}\\
\end{array}%
\right]
\in M_n(\mathcal{F}):
a^{i,j}=a^{j,i},
i,j=1,2,\dots ,n\}
$$
is a Jordan algebra with respect to the addition and the Jordan multiplication
$$
a\cdot b=\frac{1}{2}(ab+ba), a, b\in H_n(\mathcal{F}).
$$

Let $\bar{e}_{i,j}=e_{i,j}+e_{j,i}$ and $\bar{a}_{i,j}=\{e_{i,i}ae_{j,j}\}=(e_{i,i}a)e_{j,j}+e_{i,i}(ae_{j,j})$
for every $a\in H_n(\mathcal{F})$ and distinct $i$, $j$ in $\{1,2,\dots ,n\}$.

\begin{lemma}  \label{3.0}
Let $D=\sum_{k=1}^mD_{a_k,b_k}$ be an inner derivation on $H_n(\mathcal{F})$,
generated by $a_1$, $a_2$, $\dots,$ $a_m$, $b_1$, $b_2$, $\dots,$ $b_m\in H_n(\mathcal{F})$.
Then
$$
\sum_{k=1}^m [a_k,b_k]^{i,j}=-\sum_{k=1}^m [a_k,b_k]^{j,i}, i,j=1,2,\dots,n,
$$
i.e. $\sum_{k=1}^m [a_k,b_k]$ is a skew-symmetric matrix.
\end{lemma}

\begin{proof}
Indeed, let $i$, $j$ be arbitrary indices in $\{1,2,\dots ,n\}$. Then
for every $k\in\{1,2,\dots ,m\}$ we have
$$
[a_k,b_k]^{i,j}=\sum_{l=1}^n a_k^{i,l}b_k^{l,j}-\sum_{l=1}^n b_k^{i,l}a_k^{l,j},
$$
and
$$
a_k^{i,l}b_k^{l,j}=a_k^{l,i}b_k^{j,l}=b_k^{j,l}a_k^{l,i}, b_k^{i,l}a_k^{l,j}=b_k^{l,i}a_k^{j,l}=a_k^{j,l}b_k^{l,i}, l=1,2,\dots,n,
$$
since $a_k$ and $b_k$ are symmetric matrices.

Hence
$$
[a_k,b_k]^{i,j}=\sum_{l=1}^n a_k^{i,l}b_k^{l,j}-\sum_{l=1}^n b_k^{i,l}a_k^{l,j}=
$$
$$
\sum_{l=1}^n b_k^{j,l}a_k^{l,i}-\sum_{l=1}^n a_k^{j,l}b_k^{l,i}=-[a_k,b_k]^{j,i}, k=1,2,\dots,m.
$$
This completes the proof.
\end{proof}

Let $\mathcal{A}$ be a Jordan algebra.
An additive (not necessarily homogenous) map $\nabla : \mathcal{A}\to \mathcal{A}$ is called additive local inner derivation,
if for any element $x\in \mathcal{A}$ there exists an inner derivation $D$ such that
$\nabla(x)=D(x)$.

\begin{lemma}  \label{3.1}
Let $H_n(\mathcal{F})$ be the Jordan algebra
of symmetric $n\times n$ matrices over $\mathcal{F}$, $n>1$. Let $\nabla$ be an additive local inner derivation on $H_n(\mathcal{F})$.
Then for arbitrary $\lambda$, $\mu\in \mathcal{F}$ there exists an inner derivation $D$ on $H_n(\mathcal{F})$
such that
$$
\nabla (\lambda e_{i,i})=D(\lambda e_{i,i}), \nabla(\mu\bar{e}_{i,j})=D(\mu\bar{e}_{i,j}).
$$
\end{lemma}

\begin{proof}
We have
$$
\nabla(\lambda e_{i,i})+\nabla(\mu \bar{e}_{i,j})=\nabla(\lambda e_{i,i}+\mu \bar{e}_{i,j}).
$$
Let $a_1$, $a_2$, $\dots,$ $a_m$, $b_1$, $b_2$, $\dots,$ $b_m$, $c_1$, $c_2$, $\dots,$ $c_m$, $d_1$, $d_2$, $\dots,$ $d_m$ be elements in $H_n(\mathcal{F})$ such that
$$
\nabla(\lambda e_{i,i})=\sum_{k=1}^m D_{a_k,b_k}(\lambda e_{i,i}),
\nabla(\mu \bar{e}_{i,j})=\sum_{k=1}^m D_{c_k,d_k}(\mu \bar{e}_{i,j}).
$$
We have
$$
\sum_{k=1}^m D_{a_k,b_k}(\lambda e_{i,i})=D_{\frac{1}{4}\sum_{k=1}^m [a_k,b_k]}(\lambda e_{i,i}),
$$
$$
\sum_{k=1}^m D_{c_k,d_k}(\mu \bar{e}_{i,j})=D_{\frac{1}{4}\sum_{k=1}^m [c_k,d_k]}(\mu \bar{e}_{i,j}).
$$
Let $p_1$, $p_2$, $\dots,$ $p_m$, $q_1$, $q_2$, $\dots,$ $q_m$ be elements in $H_n(\mathcal{F})$ such that
$$
\nabla(\lambda e_{i,i}+\mu \bar{e}_{i,j})=\sum_{k=1}^m D_{p_k,q_k}(\lambda e_{i,i}+\mu \bar{e}_{i,j}).
$$
Then
$$
\nabla(\lambda e_{i,i}+\mu \bar{e}_{i,j})=D_{\frac{1}{4}\sum_{k=1}^m [p_k,q_k]}(\lambda e_{i,i}+\mu \bar{e}_{i,j}).
$$
Let $d_{i,i}=\frac{1}{4}\sum_{k=1}^m [a_k,b_k]$, $d_{i,j}=\frac{1}{4}\sum_{k=1}^m [c_k,d_k]$, $d=\frac{1}{4}\sum_{k=1}^m [p_k,q_k]$. Then,
since
$$
e_{i,i}[d_{i,i}\lambda e_{i,i}-\lambda e_{i,i}d_{i,i}+d_{i,j}\mu \bar{e}_{i,j}-\mu \bar{e}_{i,j}d_{i,j}]e_{i,i}=e_{i,i}[d(\lambda e_{i,i}+\mu \bar{e}_{i,j})-(\lambda e_{i,i}+\mu \bar{e}_{i,j})d]e_{i,i}
$$
we have
$$
d_{i,j}^{i,j}-d_{i,j}^{j,i}=d^{i,j}-d^{j,i},\,\,i.e \,\,d_{i,j}^{i,j}=d^{i,j}=d_{i,j}^{j,i}=d^{j,i}
$$
by lemma \ref{3.0}.
 From the equality
$$
e_{i,i}[d_{i,i}\lambda e_{i,i}-\lambda e_{i,i}d_{i,i}+d_{i,j}\mu \bar{e}_{i,j}-\mu \bar{e}_{i,j}d_{i,i}]e_{j,j}=e_{i,i}[d(\lambda e_{i,i}+\mu \bar{e}_{i,j})-(\lambda e_{i,i}+\mu \bar{e}_{i,j})d]e_{j,j}
$$
it follows that
$$
\lambda d_{i,i}^{i,j}+\mu d_{i,j}^{i,i}-\mu d_{i,j}^{j,j}=\lambda d^{i,j}+\mu d^{i,i}-\mu d^{j,j}
$$
and by lemma \ref{3.0}  we have $d_{i,j}^{i,i}=d_{i,j}^{j,j}=d^{i,i}=d^{j,j}=0$. Hence
$$
d_{i,i}^{i,j}=d^{i,j}.
$$
Therefore
$$
d_{i,i}^{i,j}=d_{i,j}^{i,j}=d_{i,i}^{j,i}=d_{i,j}^{j,i}.
$$

Let $e=e_{i,i}+e_{j,j}$. Then as in the proof of lemma \ref{1} we get
$$
(1-e)de_{i,i}=(1-e)d_{i,j}e_{i,i}.
$$
Hence $e_{i,i}d(1-e)=e_{i,i}d_{i,j}(1-e)$, since $d_{i,j}$ and $d$ are skew-symmetric matrices.

Therefore
$$
\nabla(\lambda e_{i,i})=e_{i,i}d_{i,i}\lambda e_{i,i}+e_{j,j}d_{i,i}\lambda e_{i,i}+(1-e)d_{i,i}\lambda e_{i,i}-\lambda e_{i,i}d_{i,i}e_{j,j}-
$$
$$
\lambda e_{i,i}d_{i,i}e_{i,i}-\lambda e_{i,i}d_{i,i}(1-e)=
$$
$$
e_{j,j}d_{i,j}\lambda e_{i,i}-\lambda e_{i,i}d_{i,j}e_{j,j}+e_{i,i}d_{i,j}\lambda e_{i,i}+(1-e)d_{i,j}\lambda e_{i,i}-
$$
$$
\lambda e_{i,i}d_{i,j}e_{i,i}-\lambda e_{i,i}d_{i,j}(1-e)=
$$
$$
d_{i,j}\lambda e_{i,i}-\lambda e_{i,i}d_{i,j}=D_{d_{i,j}}(\lambda e_{i,i}).
$$
since
$$
e_{i,i}d_{i,i}\lambda e_{i,i}-\lambda e_{i,i}d_{i,i}e_{i,i}=e_{i,i}d_{i,j}\lambda e_{i,i}-\lambda e_{i,i}d_{i,j}e_{i,i}=0,
$$
$$
d_{i,j}^{i,j}=d_{i,i}^{i,j}.
$$
This completes the proof.
\end{proof}

\begin{lemma} \label{3.2}
 Let $\nabla$ be an additive local inner derivation on $H_n(\mathcal{F})$.
Then for arbitrary $\lambda$, $\mu\in \mathcal{F}$ and each index $i$ there exists an inner derivation $D$ on $H_n(\mathcal{F})$
such that
$$
\nabla(\lambda e_{i,i})=D(\lambda e_{i,i}), \nabla(\mu e_{i,i})=D(\mu e_{i,i}).
$$
\end{lemma}

\begin{proof}
By Lemma \ref{3.1} for arbitrary $\nu\in \mathcal{F}$ there exist $a_1$, $a_2$, $\dots,$ $a_m$, $b_1$, $b_2$, $\dots,$ $b_m$ in $H_n(\mathcal{F})$ such that
$$
\nabla(\lambda e_{i,i})=\sum_{k=1}^m D_{a_k,b_k}(\lambda e_{i,i}),
\nabla(\nu \bar{e}_{i,j})=\sum_{k=1}^m D_{a_k,b_k}(\nu \bar{e}_{i,j}).
$$
Let $a=\frac{1}{4}\sum_{k=1}^m [a_k,b_k]$. Then
$$
\nabla(\lambda e_{i,i})=D_a(\lambda e_{i,i}), \nabla(\nu \bar{e}_{i,j})=D_a(\nu \bar{e}_{i,j}).
$$
Similarly, for arbitrary $\nu\in \mathcal{F}$ there exist $b\in M_n(\mathcal{F})$ such that
$$
\nabla(\mu e_{i,i})=D_b(\mu e_{i,i}), \nabla(\nu \bar{e}_{i,j})=D_b(\nu \bar{e}_{i,j}).
$$
We have
$$
D_a(\nu \bar{e}_{i,j})=D_b(\nu \bar{e}_{i,j}).
$$
 From
$$
e_{j,j}(D_a(\nu \bar{e}_{i,j}))=e_{j,j}(D_b(\nu \bar{e}_{i,j}))
$$
it follows that
$$
e_{j,j}ae_{i,i}=e_{j,j}be_{i,i}, e_{i,i}ae_{j,j}=e_{i,i}be_{j,j}.
$$
Let $e=e_{i,i}+e_{j,j}$. Then by the equalities
$$
(1-e)(D_a(\nu \bar{e}_{i,j}))=(1-e)(D_b(\nu \bar{e}_{i,j})), (D_a(\nu \bar{e}_{j,i}))(1-e)=(D_b(\nu \bar{e}_{j,i}))(1-e)
$$
we have $(1-e)ae_{i,i}=(1-e)be_{i,i}$ and $e_{i,i}a(1-e)=e_{i,i}b(1-e)$ respectively.

Therefore
$$
D_b(\mu e_{i,i})=b\mu e_{i,i}-\mu e_{i,i}b=
$$
$$
e_{i,i}b\mu e_{i,i}+e_{j,j}b\mu e_{i,i}+(1-e)b\mu e_{i,i}-\mu e_{i,i}be_{i,i}-\mu e_{i,i}be_{j,j}-\mu e_{i,i}b(1-e)=
$$
$$
e_{j,j}b\mu e_{i,i}-\mu e_{i,i}be_{j,j}+(1-e)b\mu e_{i,i}-\mu e_{i,i}b(1-e)=
$$
$$
e_{j,j}a\mu e_{i,i}-\mu e_{i,i}ae_{j,j}+(1-e)a\mu e_{i,i}-\mu e_{i,i}a(1-e)=
$$
$$
e_{j,j}a\mu e_{i,i}-\mu e_{i,i}ae_{j,j}+e_{i,i}a\mu e_{i,i}-\mu e_{i,i}ae_{i,i}+(1-e)a\mu e_{i,i}-\mu e_{i,i}a(1-e)=
$$
$$
a\mu e_{i,i}-\mu e_{i,i}a=D_a(\mu e_{i,i}).
$$
This completes the proof.
\end{proof}

Similarly we can prove the following lemma using lemmas \ref{4} and \ref{6}.

\begin{lemma} \label{3.3}
 Let $\nabla$ be an additive local inner derivation on $H_n(\mathcal{F})$.
Then for arbitrary $\lambda$, $\mu\in \mathcal{F}$ and each pair $i$, $j$ of distinct indices there exists an inner derivation $D$ on $H_n(\mathcal{F})$
such that
$$
\nabla(\lambda \bar{e}_{i,j})=D(\lambda \bar{e}_{i,j}), \nabla(\mu \bar{e}_{i,j})=D(\mu \bar{e}_{i,j}).
$$
\end{lemma}

 Now we have the following

\begin{lemma} \label{3.4}
 For arbitrary $\lambda$, $\mu\in \mathcal{F}$ and each pair $i$, $j$ of distinct indices we have that
$$
\nabla(\lambda e_{i,i})=\lambda \nabla(e_{i,i}), \nabla(\mu \bar{e}_{i,j})=\mu \nabla(\bar{e}_{i,j}).
$$
\end{lemma}

\begin{proof}
By lemma \ref{3.2} there exists an element $a\in M_n(\mathcal{F})$ such that
$$
\nabla(\lambda e_{i,i})=D_a(\lambda e_{i,i}), \nabla(e_{i,i})=D_a(e_{i,i}).
$$
Hence
$$
\nabla(\lambda e_{i,i})=D_a(\lambda e_{i,i})=\lambda D_a(e_{i,i})=\lambda\nabla(e_{i,i}).
$$
The proof of the second equality is similar.
\end{proof}

\begin{lemma} \label{3.5}
Let $\nabla$ be an additive local inner derivation on $H_n(\mathcal{F})$.
Then for any indices $i$, $j$, $k$, $l$, satisfying $\{i,j\}\neq \{k,l\}$, there exists an inner derivation $D$ on $H_n(\mathcal{F})$ such that
$$
\nabla(\bar{e}_{i,j})=D(\bar{e}_{i,j}), \nabla(\bar{e}_{k,l})=D(\bar{e}_{k,l}).
$$
\end{lemma}

\begin{proof}
Let $i$, $j$, $k$ be pairwise distinct indices. Then there exist inner derivation $D_1$, $D_2$ on $H_n(\mathcal{F})$ such that
$$
\nabla(\bar{e}_{i,j})=D_1(\bar{e}_{i,j}), \nabla(\bar{e}_{j,j})=D_1(\bar{e}_{j,j}),
$$
$$
\nabla(\bar{e}_{j,k})=D_2(\bar{e}_{j,k}), \nabla(\bar{e}_{j,j})=D_2(\bar{e}_{j,j}).
$$
Hence there exist elements $a$, $b\in M_n(\mathcal{F})$ such that
$$
\nabla(\bar{e}_{i,j})=D_a(\bar{e}_{i,j}), \nabla(e_{j,j})=D_a(e_{j,j}),
$$
$$
\nabla(\bar{e}_{j,k})=D_b(\bar{e}_{j,k}), \nabla(e_{j,j})=D_b(e_{j,j}).
$$
So
$$
D_a(e_{j,j})=D_a(e_{j,j})
$$
and
$$
(1-e_{j,j})ae_{j,j}=(1-e_{j,j})be_{j,j},  e_{j,j}a(1-e_{j,j})=e_{j,j}b(1-e_{j,j}).
$$
Since $e_{j,j}ae_{j,j}=e_{j,j}be_{j,j}=0$ we have
$$
ae_{j,j}=be_{j,j},  e_{j,j}a=e_{j,j}b, e_{i,i}ae_{j,j}=e_{i,i}be_{j,j},  e_{j,j}ae_{i,i}=e_{j,j}be_{i,i}.
$$
Also we have
$$
\nabla(\bar{e}_{i,j})+\nabla(\bar{e}_{j,k})=\nabla(\bar{e}_{i,j}+\bar{e}_{j,k})
$$
Let $d\in M_n(\mathcal{F})$ be an element such that
$$
\nabla(\bar{e}_{i,j}+\bar{e}_{j,k})=D_d(\bar{e}_{i,j}+\bar{e}_{j,k}).
$$
Then
$$
D_a(\bar{e}_{i,j})+D_b(\bar{e}_{j,k})=D_d(\bar{e}_{i,j}+\bar{e}_{j,k}).
$$
From
$$
e_{i,i}[D_a(\bar{e}_{i,j})+D_b(\bar{e}_{j,k})]e_{j,j}=e_{i,i}[D_d(\bar{e}_{i,j}+\bar{e}_{j,k})]e_{j,j}
$$
and
$$
e_{j,j}[D_a(\bar{e}_{i,j})+D_b(\bar{e}_{j,k})]e_{k,k}=e_{j,j}[D_d(\bar{e}_{i,j}+\bar{e}_{j,k})]e_{k,k}
$$
it follows that $e_{i,i}be_{k,k}=e_{i,i}de_{k,k}$ and $e_{i,i}ae_{k,k}=e_{i,i}de_{k,k}$ respectively.
Hence $e_{i,i}ae_{k,k}=e_{i,i}be_{k,k}$ and $e_{k,k}ae_{i,i}=e_{k,k}be_{i,i}$.
Now, let $e=e_{i,i}+e_{j,j}$. Then
$$
(1-e-e_{k,k})be_{i,i}\bar{e}_{j,k}=\bar{e}_{j,k}(1-e-e_{k,k})be_{i,i}=0,
$$
$$
e_{i,i}b(1-e-e_{k,k})\bar{e}_{j,k}=\bar{e}_{j,k}e_{i,i}b(1-e-e_{k,k})=0.
$$
So we may assume that
$$
(1-e-e_{k,k})be_{i,i}=(1-e-e_{k,k})ae_{i,i}, e_{i,i}b(1-e-e_{k,k})=e_{i,i}a(1-e-e_{k,k}).
$$
Hence
$$
(1-e)be_{i,i}=(1-e)ae_{i,i}, e_{i,i}b(1-e)=e_{i,i}a(1-e)
$$
since $e_{i,i}ae_{k,k}=e_{i,i}be_{k,k}$ and $e_{k,k}ae_{i,i}=e_{k,k}be_{i,i}$. Therefore
$$
\nabla(\bar{e}_{i,j})=D_a(\bar{e}_{i,j})=a\bar{e}_{i,j}-\bar{e}_{i,j}a=
$$
$$
ae_{j,i}+(1-e)ae_{i,j}+e_{i,i}ae_{i,j}+e_{j,j}ae_{i,j}-
$$
$$
e_{i,j}a-e_{j,i}a(1-e)-e_{j,i}ae_{i,i}-e_{j,i}ae_{j,j}=
$$
$$
be_{j,i}+(1-e)be_{i,j}+e_{i,i}be_{i,j}+e_{j,j}be_{i,j}-
$$
$$
e_{i,j}b-e_{j,i}b(1-e)-e_{j,i}be_{i,i}-e_{j,i}be_{j,j}=D_b(\bar{e}_{i,j}).
$$

Let $i$, $j$, $k$, $l$ be pairwise distinct indices and let
$a$, $b$ be elements in $M_n(\mathcal{F})$ such that
$$
\nabla (\bar{e}_{i,j})=D_{a}(\bar{e}_{i,j}), \nabla (\bar{e}_{k,l})=D_{b}(\bar{e}_{k,l})
$$
and let $e=e_{i,i}+e_{j,j}$, $f=e_{k,k}+e_{l,l}$.
Then we have
$$
(1-f)b(1-f)\bar{e}_{k,l}=\bar{e}_{k,l}(1-f)b(1-f)=0.
$$
So we may take $(1-f)b(1-f)=(1-f)a(1-f)$. In particular,
$$
e_{i,i}b(1-f)=e_{i,i}a(1-f), (1-f)be_{i,i}=(1-f)ae_{i,i},
$$
$$
e_{j,j}b(1-f)=e_{j,j}a(1-f), (1-f)be_{j,j}=(1-f)ae_{j,j}.
$$
Then $D_d(\bar{e}_{i,j}+\bar{e}_{k,l})=D_{a}(\bar{e}_{i,j})+D_{b}(\bar{e}_{k,l})$ and by the equalities
$$
e_{i,i}(D_d(\bar{e}_{i,j}+\bar{e}_{k,l}))e_{k,k}=e_{i,i}(D_{a}(\bar{e}_{i,j})+D_{b}(\bar{e}_{k,l}))e_{k,k},
$$
$$
e_{j,j}(D_d(\bar{e}_{i,j}+\bar{e}_{k,l}))e_{l,l}=e_{j,j}(D_{a}(\bar{e}_{i,j})+D_{b}(\bar{e}_{k,l}))e_{l,l}
$$
we have
$$
d^{i,l}-d^{j,k}=b^{i,l}-a^{j,k}, d^{j,k}-d^{i,l}=b^{j,k}-a^{i,l},
$$
respectively. Hence $b^{i,l}-a^{j,k}=a^{i,l}-b^{j,k}$.
Also, by the equalities
$$
\nabla(\bar{e}_{i,j}-\bar{e}_{k,l})=\nabla(\bar{e}_{i,j})-\nabla(\bar{e}_{k,l}),
$$
$$
e_{i,i}(D_d(\bar{e}_{i,j}-\bar{e}_{k,l}))e_{k,k}=e_{i,i}(D_{a}(\bar{e}_{i,j})-D_{b}(\bar{e}_{k,l}))e_{k,k},
$$
$$
e_{j,j}(D_d(\bar{e}_{i,j}-\bar{e}_{k,l}))e_{l,l}=e_{j,j}(D_{a}(\bar{e}_{i,j})-D_{b}(\bar{e}_{k,l}))e_{l,l}
$$
we have
$$
d^{i,l}+d^{j,k}=b^{i,l}+a^{j,k}, d^{j,k}+d^{i,l}=b^{j,k}+a^{i,l},
$$
respectively. Hence $b^{i,l}+a^{j,k}=a^{i,l}+b^{j,k}$. Therefore
$$
b^{i,l}=a^{i,l}, a^{j,k}=b^{j,k}\,\,and\,\,b^{l,i}=a^{l,i}, a^{k,j}=b^{k,j}.
$$
Similarly we get
$$
b^{i,k}=a^{i,k}, a^{j,l}=b^{j,l}\,\,and\,\,b^{k,i}=a^{k,i}, a^{l,j}=b^{l,j}.
$$
Therefore
$$
\nabla(\bar{e}_{i,j})=e_{k,k}a\bar{e}_{i,j}+e_{l,l}a\bar{e}_{i,j}+(1-f)a\bar{e}_{i,j}-\bar{e}_{i,j}ae_{l,l}-
$$
$$
\bar{e}_{i,j}ae_{k,k}-\bar{e}_{i,j}a(1-f)=
$$
$$
e_{k,k}b\bar{e}_{i,j}+e_{l,l}b\bar{e}_{i,j}+(1-f)b\bar{e}_{i,j}-\bar{e}_{i,j}be_{l,l}-
$$
$$
\bar{e}_{i,j}be_{k,k}-\bar{e}_{i,j}b(1-f)=
$$
$$
b\bar{e}_{i,j}-\bar{e}_{i,j}b=D_{b}(\bar{e}_{i,j}).
$$
This completes the proof.
\end{proof}

\begin{lemma} \label{3.6}
Let $\nabla :H_n(\mathcal{F})\to H_n(\mathcal{F})$ be an additive local inner derivation.
Then there exists a skew-symmetric matrix $a\in M_n(\mathcal{F})$ such that for any indices $i$, $j$ we have
$$
\nabla (\bar{e}_{i,j})=D_a(\bar{e}_{i,j}).
$$
\end{lemma}

\begin{proof}
By the previous lemmas we can repeat the proof of lemma 14 in \cite{AA3}
and get the statement of the above lemma. The proof is complete.
\end{proof}

\begin{theorem} \label{3.7}
Let  $H_n(\mathcal{F})$ be the Jordan algebra
of $n\times n$ symmetric matrices over $\mathcal{F}$, $n>1$. Then any additive local inner derivation on
$H_n(\mathcal{F})$ is a derivation.
\end{theorem}

\begin{proof}
Let $\nabla :H_n(\mathcal{F})\to H_n(\mathcal{F})$ be an additive local inner derivation. Then
by lemma \ref{3.6} there exists a skew-symmetric matrix $a\in M_n(\mathcal{F})$ such that for any indices $i$, $j$ we have
$$
\nabla (\bar{e}_{i,j})=D_a(\bar{e}_{i,j}).
$$
Let $x$ be an arbitrary element in $H_n(\mathcal{F})$.
Then $x=\sum_{k,l=1}^n x^{k,l}e_{k,l}$ and by lemma \ref{3.4}
$$
\nabla(x)=\sum_{k,l=1, k\leq l}^n \nabla(x^{k,l}\bar{e}_{k,l})=\sum_{k,l=1, k\leq l}^n x^{k,l}\nabla(\bar{e}_{k,l})=
$$
$$
\sum_{k,l=1, k\leq l}^n x^{k,l}D_a(\bar{e}_{k,l})=D_a(x).
$$
Hence $\nabla$ is a derivation on $H_n(\mathcal{F})$ (to be more precise, $\nabla$ is a spatial derivation implemented by a skew-symmetric matrix $a\in M_n(\mathcal{F})$) . The proof is complete.
\end{proof}

Let $\mathcal{J}$ be a Jordan ring. An additive map $\nabla : \mathcal{J}\to \mathcal{J}$ is called local derivation, if for
any element $x\in \mathcal{J}$ there exists a derivation
$D:\mathcal{J}\to \mathcal{J}$ such that $\nabla (x)=D(x)$.

An additive map $\nabla : \mathcal{J}\to \mathcal{J}$ is called local inner derivation,
if for any element $x\in \mathcal{J}$ there exists an inner derivation $D$ on the Jordan ring $\mathcal{J}$ such that
$\nabla(x)=D(x)$.

By the proofs of the previous lemmas of the present section we have the following lemma.

\begin{lemma} \label{3.8}
Let $\Re$ be a finite ring generated by the identity element
or the ring of integers, $H_n(\Re)$ be the Jordan ring
of $n\times n$ symmetric matrices over $\Re$, $n>1$, and let $\nabla :H_n(\Re)\to H_n(\Re)$ be a local inner derivation.
Then there exists a skew-symmetric matrix $a\in M_n(\lambda\Re)$
such that for any indices $i$, $j$ we have
$$
\nabla(\bar{e}_{i,j})=D_a(\bar{e}_{i,j}).
$$
\end{lemma}

By lemma \ref{3.8} and by the definition of a local inner derivation we have the following theorem.

\begin{theorem} \label{3.9}
Let $\Re$ be a finite ring generated by the identity element
or the ring of integers, $H_n(\Re)$ be the Jordan ring
of $n\times n$ symmetric matrices over $\Re$, $n>1$. Then for every local inner derivation $\nabla$ on
the Jordan ring $H_n(\Re)$ there exists a skew-symmetric matrix $a\in M_n(\Re)$ such that
$$
\nabla(x)=D_a(x), x\in H_n(\Re),
$$
i.e. $\nabla$ is a derivation on
the Jordan ring $H_n(\Re)$.
\end{theorem}

\end{document}